\documentclass[12pt]{amsart}

\usepackage{amsmath, amsthm, amssymb}
\usepackage{geometry}
\usepackage[
    colorlinks=true,
    allcolors=blue,
    unicode=true
]{hyperref}
\usepackage{graphicx}
\usepackage{caption}
\usepackage{subcaption}
\usepackage{mathtools}
\usepackage[capitalise]{cleveref}
\usepackage[style=british]{csquotes}
\usepackage{parskip}
\usepackage{float}

\usepackage[
  backend=biber,
  style=numeric-comp,
]{biblatex}
\AtEveryBibitem{%
  \clearfield{note}%
  \clearfield{urlyear}%
}

\DeclareSourcemap{
  \maps[datatype=bibtex]{
    \map{
      \step[fieldsource=doi, final]
      \step[fieldset=url, null]
    }
  }
}

\newtheorem{theorem}{Theorem}
\newtheorem{prop}{Proposition}
\newtheorem*{prop*}{Proposition}
\newtheorem*{theorem*}{Theorem}
\newtheorem{lemma}{Lemma}
\newtheorem{corollary}{Corollary}

\newcommand{\CC}{\mathbb{C}}
\newcommand{\RR}{\mathbb{R}}
\newcommand{\QQ}{\mathbb{Q}}
\newcommand{\ZZ}{\mathbb{Z}}
\newcommand{\NN}{\mathbb{N}}
\newcommand{\FF}{\mathbb{F}}

\newcommand{\acrit}{a_\text{crit}}

\DeclarePairedDelimiter{\paren}{\lparen}{\rparen}
\DeclarePairedDelimiter{\braces}{\lbrace}{\rbrace}
\DeclarePairedDelimiter{\ceil}{\lceil}{\rceil}
\DeclarePairedDelimiter{\abs}{\lvert}{\rvert}

\DeclareMathOperator{\wam}{\mathfrak{W}}
\DeclareMathOperator{\rad}{rad}

\usepackage[foot]{amsaddr}
\title[A remark on weighted average multiplicities in prime factorisation]{A remark on weighted average multiplicities \\ in prime factorisation}
\author{Viktor Mirjani\'c, Daattavya Aggarwal, Challenger Mishra}
\address{Department of Computer Science and Technology, University of Cambridge, Cambridge, United Kingdom}
\email{ \{vvm22,da579,cm2099\}@cam.ac.uk}

\begin{document}

\begin{abstract}
We study a generalisation of the quality of an ABC triple that we call the weighted average multiplicity (WAM), in which the logarithmic heights of prime factors are raised to a complex exponent $s$. The WAM is connected to the standard ABC conjecture at $s=1$. 
We show that for real part of $s$ less than $1$, WAM is unbounded over ABC triples both for integers and polynomials. For real part greater than $1$, we characterise a boundary beyond which WAM is holomorphic and bounded. In this region, we show that WAM is related to the multiplicity of the largest prime factor of the triple, a quantity that we connect with the original ABC conjecture and whose distribution we explore computationally.
\end{abstract} 
\maketitle

\section{Weighted Average Multiplities}

Given an integer $n$ with prime factorisation 
\begin{equation*}
    n=\pm \prod_{i=1}^m p_k^{e_k},
\end{equation*}
the {\em weighted average multiplicity  (WAM)} 
of its set of prime factors is defined by
\begin{equation}
\label{eq:original_wam_definition}
    \wam(n):=\frac{\sum_{k=1}^m e_k \log p_k}{\sum_{k=1}^m\log p_k}.
\end{equation}
This definition is due to Minhyong Kim. Here as in the following, the indexing of primes in a factorisation like this will be in increasing order. We will use the standard terminology {\em height} of an integer $n$ for $\log |n|$, which is roughly its information-theoretic complexity.

This note argues that this definition provides a natural way to formulate and investigate a number of  properties of factorisation, especially the difficult ideas surrounding the {\em ABC conjecture}\footnote{A proof of this conjecture was announced by Shinichi Mochizuki in 2012 in a series of four preprints, which were published in 2021 \cite{mochizuki_inter-universal_2021, mochizuki_inter-universal_2021-1, mochizuki_inter-universal_2021-2, mochizuki_inter-universal_2021-3}. In 2018, Peter Scholze and Jakob Stix visited Kyoto for discussions with Mochizuki and subsequently wrote a report \cite{scholze_why_2018} explaining what they viewed as a gap. Recently, Kirti Joshi released a series of preprints \cite{joshi_construction_2023, joshi_construction_2025, joshi_construction_2025-1, joshi_construction_2025-2, joshi_construction_2025-3, joshi_construction_2025-4, joshi_final_2025} claiming to have fixed the gap. The discussion in our paper is independent of the status of this long disagreement, other than the fact that we refer to the key inequality as a conjecture.}. The reader will notice right away that
\begin{equation}
\label{eq:original_wam_via_log_rad}
    \wam(n)=\frac{\log |n|}{\log \rad |n|},
\end{equation}
where
\begin{equation*}
    \rad |n|:=\prod_{k=1}^m p_k,
\end{equation*}
is the radical of $n$, or the number obtained from $n$ by removing all multiplicities. This is the connection to ABC as usually formulated. We call an {\em ABC triple} a triple of coprime integers $a,b,c$  such that $a+b=c$.
It is conjectured that for any $\epsilon>0$, there are only finitely many ABC triples such that
\begin{equation}
\label{eq:wam_abc_conjectura}
    \wam(abc)>3+\epsilon.
\end{equation}
We note that this is implied by the standard form of the ABC conjecture, which is sometimes expressed in terms of the `quality' $q$ of the triple:
\begin{equation}
q(a,b,c):=\frac{\log \max(|a|,|b|, |c|)}{\log \rad |abc|},
\end{equation}
conjectured to be less than $1+\epsilon$ for all but finitely many triples.
Our statement  might be strictly weaker than the usual conjecture and, in particular, it is not at all clear that the bound $3+\epsilon$ is optimal.  Nonetheless, since the precise exponent is of little conceptual importance, we will stick to this formulation and confine our  investigation to  the WAM rather than the quality in this paper.
In any case, it appears natural that multiplicities are weighted by the heights of the primes they contribute, and the ABC conjecture suggests that this weighting is what keeps the average multiplicities of ABC triples bounded: If small primes occur in $abc$ with large multiplicities, compensation will come from large primes with small multiplicities. The meaning of such informal accounts of the conjecture are made precise by the expression \ref{eq:original_wam_definition} for the WAM, which is not as apparent with the simpler expression \ref{eq:original_wam_via_log_rad}.

Another advantage of making the weights explicit is that we can then vary them in a family to define the  function
\begin{equation}
\label{eq:wam_def}
    \wam(n,s):=\frac{\sum_{k=1}^m e_k (\log p_k)^s}{\sum_{k=1}^m(\log p_k)^s},
\end{equation}
giving the heights more or less prominence depending on the size of $s$. We allow $s$ to take on complex values, whereby we get a meromorphic function on $\CC$. A feature of this family is that while the value at $s=1$ is the original WAM, the value at zero 
\begin{equation*}
    \wam(n,0)=\frac{\sum_{k=1}^me_k}{m}
\end{equation*}
is the normal average multiplicity, and
\begin{equation*}
    \lim_{\Re(s) \to +\infty} \wam(n,s)=e_m,    
\end{equation*}
the multiplicity of the largest prime factor.

Note the inequality
\begin{equation*}
    \wam(n,1)=\frac{\sum_{k=1}^me_k(\log p_k)}{\sum_{k=1}^m(\log p_k)}\geq \frac{e_m\log p_m}{\sum_{k=1}^m(\log p_k)}\geq \frac{e_m}{m}
\end{equation*}
giving us the crude bound
\begin{equation*}
    e_m\leq \wam(n,1)\omega(n).
\end{equation*}
on $e_m$. Here, $\omega(n)$ is the usual notation for the number of distinct prime divisors of $n$.
A theorem of Erdos and Kac  \cite{ErdosKac} implies that $e_m$ will tend not to be much larger than $\wam(n,1)$. It is not clear how large $e_m$ can be for ABC triples, but it is quite easy to see that the average multiplicity $\wam(abc,0)$ is unbounded.

The definite  mathematical results of this paper are quite modest and amount to the observation that the ABC inequalities appear to be quite finely balanced.
\begin{theorem}
\label{thm:wam_divergence_integers}
If $\Re(s)<1$, then the quantity
\begin{equation*}
    \wam(abc,s)
\end{equation*}
is unbounded as we run over ABC triples of integers.
\end{theorem}
In short, we cannot hope for a bound if the weights are even slightly less than the heights of prime divisors. In fact, the ABC triples $(2^n, 1-2^n, 1)$ for large $n$ suffice to prove the theorem.

Recall the standard analogy between integers and  polynomials over finite fields. That is, let $\FF_q$ be a finite field and  let $f\in \FF_q[x]$ have factorisation
\begin{equation*}
    f(x) =\prod_{i=1}^m p_k^{e_k}
\end{equation*}
into irreducible polynomials $p_k$. We define
\begin{equation*}
    \wam(f,s):=\frac{\sum_{k=1}^m e_k (\deg p_k)^s}{\sum_{k=1}^m(\deg p_k)^s}
\end{equation*}
An ABC triple in $\FF_q[x]$ are coprime polynomials $a,b,c$ such that $c=a+b$ and at least one of $a,b,c$ have a non-zero formal derivative. In this case, it is shown by \textcite[Lemma 2]{mason_diophantine_1984} and Stothers  that
\begin{equation*}
    \wam(abc)=\wam(abc, 1)\leq 3,
\end{equation*}
so that an \textquote{ABC inequality} is not conjectural. Nonetheless, we show
\begin{theorem}
\label{thm:wam_divergence_polynomials}
If $\Re(s)<1$, then the quantity
$$\wam(abc,s)$$
is unbounded as we run over ABC triples of polynomials.
\end{theorem}
The proof here is somewhat different from the integer case and utilises an interesting counting argument. In principle, it should be possible to use the explicit ABC triples $(x^{\ell^n}, 1-x^{\ell^n}, 1)$ for a suitable  prime $\ell$. However, for the irreducibility of cyclotomic polynomials, we need to   find an odd $\ell$ such that $q$ is a primitive root modulo $\ell^2$. Such an  $\ell$ almost certainly always exists. For example, for $\FF_2[x]$, we can use $\ell=3$. However, it seems hard to prove the existence in general. 

Having asserted the finely tuned nature of the conjecture (in $\ZZ$) and the theorem (in $\FF_q[x]$), we need to ask what happens when $\Re(s)>1$, that is, when the heights are  given even more importance.

This question appears to depend delicately on the possible pole structures of the WAM coming from the denominators
\begin{equation*}
    \sum_{k=1}^m (\log p_k)^s.
\end{equation*}

In general, numerical experiments indicate that the poles are quite chaotically distributed as we run over ABC triples. However, we can still make some remarks.
We tackle this by first noticing that for any $n$, we can find $\acrit = \Re(s)$ such that $\wam(abc, s)$ is bounded in the region $\Re(s) > \acrit$. Furthermore, we show that our choice of $\acrit$ is optimal.

The structure of the paper is as follows. In \cref{sec:wam_divergence_integers,sec:wam_divergence_poly} we prove the unboundedness of $\wam(n,s)$ for both integer and polynomial triples when $\Re(s)<1$. Then, in \cref{sec:wam_analysis_numerical} we perform multiple numerical experiments with a dataset of \textcite{abcDataset} in order to explore the relationship between $\wam$ and ABC triples in the region $\Re(s)>1$. We then show the boundedness of $\wam$ for $\Re(s)>\acrit$, and present some plots to show the behaviour of the poles in the region $1<\Re(s)<\acrit$.

Looking back, we remark on the significance of computation for our understanding of $\wam$ even in the region where $\Re(s)<1$. We believe that further computation will help with questions raised by our observations on the distribution of poles of $\wam$ when $\Re(s)<\acrit$, or on the multiplicity of the largest prime factor $e_m$ over ABC triples.

\subsection{Acknowledgements}

We would like to thank Minhyong Kim for proposing the definition of weighted average multiplicity $\wam(n,s)$, and for many helpful discussions and comments. 
\section[Proof of divergence theorem for integers]{Proof of \cref{thm:wam_divergence_integers}}
\label{sec:wam_divergence_integers}

In this section we will prove \cref{thm:wam_divergence_integers}. We will show divergence of $\wam$ for integers by computing $\wam$ over $abc$ triples $(1, 2^n-1, 2^n)$ and showing that it will be unbounded as $n$ tends to infinity. Firstly, we will prove the simpler case for real $s$, and then we will generalise to complex $s$.

\subsection[Divergence for real s<1]{Divergence for real $s<1$}

Behaviour of $\wam\paren[\big]{2^n(2^n-1), s}$ depends on prime factorisation of $2^n-1$. While this factorisation is inaccessible, it will be sufficiently regular that we can place bounds on it. Let
\begin{equation*}
    2^n(2^n-1) = 2^n p_1^{e_1}\dots p_k^{e_k}
\end{equation*}
be the prime factorisation of $2^n(2^n-1)$, and $k=\omega(2^n-1)$ be the number of prime factors of $2^k-1$. Then, we compute
\begin{equation*}
    \wam\paren[\big]{2^n(2^n-1), s} = \frac{n (\log 2)^s + \sum e_i (\log p_i)^s}{(\log 2)^s + \sum (\log p_i)^s}.
\end{equation*}
By placing suitable lower and upper bounds on the numerator and denominator respectively, we will show that the numerator of this expression asymptotically grows strictly faster than the denominator, which will make $\wam$ diverge in the limit.

For the numerator we use the trivial bound
\begin{equation}
\label{eq:wam_num_real}
    n (\log 2)^s + \sum e_i (\log p_i)^s > n (\log 2)^s.
\end{equation}
Next, we will show that the denominator grows strictly slower than that. We will consider cases $s\le 0$ and $0<s<1$ separately.

\textbf{Denominator bound for $s\le 0$:}

        When $s \le 0$ we have that $(\log 2)^s > (\log p_i)^s$ for all primes $p_i$ since $p_i>2$. Therefore, the denominator satisfies 
        \begin{equation*}
            (\log 2)^s + \sum (\log p_i)^s < (k + 1) (\log 2)^s
        \end{equation*}
        From this and \cref{eq:wam_num_real} we conclude 
        \begin{equation}
        \label{eq:wam_num_denom_bound_s>0}
            \wam\paren[\big]{2^n(2^n-1), s} > \frac{n}{k+1}
        \end{equation}
        Now it suffices to show that $k$ grows much slower than $n$. Recall that $k=\omega(2^n-1)$. From prime number theorem we have the asymptotic bound that for all $\epsilon>0$ and $a > a_0$
        \begin{equation*}
            \omega(a) < (1+\epsilon)\frac{\log a}{\log \log a},
        \end{equation*}
        and even stronger bounds exist in \cite{robin_estimation_1983}. Substituting $2^n-1$ we conclude that asymptotically 
        \begin{equation*}
            k < c \frac{n}{\log n},
        \end{equation*}
        and applying this to \cref{eq:wam_num_denom_bound_s>0} shows that $\wam$ grows at least as fast as $\log n$. Note that $\omega(a) = o(\log a)$ is sufficient for the proof, but the stronger bound allows us to comment on the speed of the divergence.
        
\textbf{Denominator bound for $0 < s < 1$:}
    
        When $s < 0 < 1$, the function $x\mapsto x^s$ is concave, and we will be using Jensen's inequality 
        \begin{equation*}
             \sum \lambda_i f(x_i) \le f\paren*{\sum \lambda_i x_i}
        \end{equation*}
        for constants $\lambda_i$ satisfying $\lambda_i>0$, and $\sum \lambda_i = 1$. With this, we bound the denominator as follows
        {\allowdisplaybreaks
        \begin{align*}
             (\log 2)^s + \sum^k (\log p_i)^s &= (\log 2)^s + k \sum \frac{1}{k} (\log p_i)^s && \\
            &\le (\log 2)^s + k \paren*{\sum \frac{1}{k} \log p_i}^s && \text{Jensen's inequality for } (-)^{s} \\
            &= (\log 2)^s + k^{1-s} \paren*{\sum \log p_i}^s && \text{factor out } k \\
            &= (\log 2)^s + k^{1-s} \paren*{\log\rad (2^n-1)}^s  && \text{simplify sum}\\
            &\leq (\log 2)^s + k^{1-s} \paren*{\log(2^n-1)}^s  && \rad(x)\leq x\\
            &< (\log 2)^s + (1+\epsilon) \paren*{\frac{n}{\log n}}^{1-s} \paren*{\log(2^n-1)}^s  && \text{bound } k\\
            &< (\log 2)^s + (1+\epsilon) \paren*{\frac{n}{\log n}}^{1-s} \paren*{n\log 2}^s  && \text{by } \log(2^n-1)<\log(2^n) \\
            &< (\log 2)^s + (1+\epsilon) \frac{n \paren{\log 2}^s}{\paren{\log n}^{1-s}}  && \text{simplify} \\
            &\sim (\log 2)^s \frac{n}{\paren{\log n}^{1-s}}  && \text{asymptotic bound} \\
        \end{align*}
        }
        From this, we have that asymptotically $\wam$ is lower bound by
        \begin{equation*}
            \wam\paren[\big]{2^n(2^n-1), s} > 
            \frac{n (\log 2)^s}{(\log 2)^s \frac{n}{\paren{\log n}^{1-s}}} = 
            (\log n)^{1-s},
        \end{equation*}
        which diverges since $s<1$.

This proof is formalized in Lean 4 \cite{moura_lean_2021}, and can be found at \url{https://github.com/mirjanic/wam_lean}.

\subsection[Divergence for complex s with real part <1]{Divergence for complex $s$ with $\Re(s)<1$}

If $s$ is complex then so is $\wam$, so we will prove that it diverges in magnitude. Specifically, we will show that 
\begin{equation}
\label{eq:complex_bound_goal}
    \abs*{ \wam \paren[\big]{2^n(2^n-1), s} } > \frac{1}{2} \paren*{ 1 - \paren*{\frac{\log 2}{\log 3} }^{1 - \Re(s)} } \wam \paren[\Big]{ 2^n (2^n-1), \Re(s) }.
\end{equation}
Note that the given constant factor is positive. Together with the result from the previous section, this implies that $\wam$ does diverge in magnitude for $\Re(s)<1$. To show \cref{eq:complex_bound_goal} we begin with the following lemma.
\begin{lemma}
\label{lemma:pf_complex_bound}
    The prime factorisation of $2^n-1$ satisfies
    \begin{equation*}
        \sum e_i \abs*{\paren*{\log p_i}^s}< n \paren{\log 3}^{\Re(s)-1} \log 2 
    \end{equation*}
    
    \begin{proof}
        We consider two cases depending on whether $\Re(s)>0$ or not. 

        \textbf{Case $\Re(s) \le 0$:}

                Since $p_i \ge 3$ then $\log p_i \ge \log 3$. But, because $\Re(s)<0$ we have $(\log p_i)^{\Re(s)} \le (\log 3)^{\Re(s)}$. From this we conclude
                {\allowdisplaybreaks
                \begin{align*}
                    \sum e_i \abs{ (\log p_i)^s } 
                    &= \sum e_i (\log p_i)^{\Re(s)} && \abs{(-)^s} \mapsto (-)^{\Re(s)}\\
                    &\le \sum e_i (\log 3)^{\Re(s)} && \text{by above} \\
                    &= (\log 3)^{\Re(s)} \Omega(2^n-1) && \text{by } \sum e_i = \Omega(2^n-1) \\
                    &\le (\log 3)^{\Re(s)} \frac{\log(2^n-1)}{\log 3} &&  \text{by } \Omega(x) \le \frac{\log x}{\log p_\text{min}(x)} \\
                    &< (\log 3)^{\Re(s)} \frac{n \log 2}{\log 3} && \text{by } \log(2^n-1) < \log(2^n) \\
                    &= n \paren{\log 3}^{\Re(s)-1} \log 2 
                \end{align*}
                }
            \textbf{Case $0 < \Re(s) < 1$:}

                In this case the function $x\mapsto x^s$ is concave, and we will use Jensen's inequality. Let $K = \Omega(2^n-1)$.
                {\allowdisplaybreaks
                \begin{align*}
                    \sum e_i \abs{ (\log p_i)^s } 
                    &= K \sum \frac{e_i}{K} (\log p_i)^{\Re(s)} && \abs{(-)^s} \mapsto (-)^{\Re(s)}\\
                    &< K \paren*{\sum \frac{e_i}{K} \log p_i}^{\Re (s)} && \text{Jensen's inequality for } (-)^{\Re(s)}  \\
                    &= K^{1-\Re(s)} \paren*{\sum e_i \log p_i}^{\Re (s)} && \text{factor out } K \\
                    &= K^{1-\Re(s)} \paren*{\log(2^n-1)}^{\Re (s)} && \text{simplify sum} \\
                    &\le \paren*{\frac{\log (2^n-1)}{\log 3}}^{1-\Re(s)} \paren*{\log(2^n-1)}^{\Re (s)}  &&  \text{by } \Omega(x) \le \frac{\log x}{\log p_\text{min}(x)} \\
                    &= (\log 3)^{\Re(s)-1} \log(2^n-1)  &&  \text{factor} \\
                    &< n (\log 3)^{\Re(s)-1} \log 2 &&  \text{by } \log(2^n-1) < \log(2^n) 
                \end{align*}
                }
    \end{proof}
\end{lemma}
In particular, we note that this implies
\begin{equation}
\label{eq:pf_complex_bound}
    \sum e_i \abs*{\paren*{\log p_i}^s} = \sum e_i (\log p_i)^{\Re(s)} < n (\log 2)^{\Re(s)} = n \abs{(\log 2)^s}    
\end{equation}
Next, we apply triangle inequalities to the numerator and denominator of $\abs{\wam}$:
{\allowdisplaybreaks
\begin{align*}            
    \abs*{\wam\paren[\big]{2^n(2^n-1),s}}&= 
    \abs*{\frac{ n (\log 2)^s +\sum e_i(\log p_i)^s} {(\log 2)^s + \sum (\log p_i)^s}} && \\
    &\ge \frac{n \abs{(\log 2)^s}-\sum e_i\abs*{(\log p_i)^s}}{\abs{(\log 2)^s} + \sum \abs*{(\log p_i)^s}} && \text{triangle inequalities} \\
    &> \frac{n \abs{(\log 2)^s}- \paren{\log 3}^{\Re(s)-1} \log(2^n-1) }{\abs{(\log 2)^s} + \sum \abs*{(\log p_i)^s}} && \text{by \cref{lemma:pf_complex_bound}} \\
    &> \frac{n \abs{(\log 2)^s}- n(\log 3)^{\Re(s)-1} \log 2 }{\abs{(\log 2)^s} + \sum \abs*{(\log p_i)^s}} && \text{by } \log(2^n-1)<\log(2^n) \\
    &= \frac{n (\log 2)^{\Re(s)} - n(\log 3)^{\Re(s)-1} \log 2 }{(\log 2)^{\Re(s)} + \sum (\log p_i)^{\Re(s)}} && \abs{(-)^s} \mapsto (-)^{\Re(s)} \\
    &=  \paren*{ 1 - \paren*{\frac{\log 2}{\log 3} }^{1 - \Re(s)} } \frac{n (\log 2)^{\Re(s)} }{(\log 2)^{\Re(s)} + \sum (\log p_i)^{\Re(s)}} && \text{simplify numerator} \\
    &>  \frac{1}{2} \paren*{ 1 - \paren*{\frac{\log 2}{\log 3} }^{1 - \Re(s)} }\wam \paren[\Big]{ 2^n (2^n-1), \Re(s) } && \text{by \cref{eq:pf_complex_bound}} \\
\end{align*}
}
Since the latter diverges, we are done.
\section[Proof of divergence theorem for polynomials]{Proof of \cref{thm:wam_divergence_polynomials}}
\label{sec:wam_divergence_poly}

In this section we will prove \cref{thm:wam_divergence_polynomials}, which states that $\wam$ diverges over $abc$ triples of polynomials for $\Re(s)<1$. We will first prove the theorem for $\QQ[x]$, and then for $\FF_p[x]$. These proofs will be significantly less demanding compared to the integer case. In both cases, we will construct triples that have few irreducible factors but large multiplicities, and this will make $\wam$ diverge in the appropriate limit.

\subsection{Proof for \texorpdfstring{$\QQ[x]$}{rational polynomials}}

Consider the following polynomial triple: $(P,Q,R) = (1, x^p-1,x^p)$, for some prime $p$. This is a valid $abc$ triple since clearly $P+Q=R$, and they are not all constant. We will proceed to compute $\wam$ for these triples and show that $\lim_{p\to+\infty} \wam(PQR,s)$ diverges for $\Re(s)<1$.

To compute $\wam$ we need to factorise $PQR$, and this amounts to factorising $Q$. We have 
\begin{equation*}
    Q(x) = (x - 1) (1 + x + \dots + x^{p-1}) = (x - 1) \Phi_p (x),
\end{equation*}
where $\Phi_p$ is the cyclotomic polynomial, which is known to be irreducible in $\QQ[x]$. Thus, we have $PQR=(x-1)x^p \Phi_p(x)$, so
\begin{equation*}
    \wam(PQR,s)=
    \frac{\sum_{k=1}^m e_k (\deg p_k)^s}{\sum_{k=1}^m(\deg p_k)^s} = 
    \frac{\overbrace{1\times1^s}^\text{from $x-1$} + \overbrace{1\times p^s}^\text{from $\Phi_p(x)$} + \overbrace{p\times1^s}^\text{from $x^p$}}{1^s + p^s + 1^s} = 
    \frac{p+p^s+1}{p^s+2}.
\end{equation*}
When $s\in\RR$ and $s<1$, this expression diverges since $p \gg p^s$ in the limit of $p\to+\infty$ over primes $p$. For complex $s$ with $\Re(s)<1$, the same argument holds except we consider the magnitude of $\wam$ instead.

\subsection{Proof for \texorpdfstring{$\FF_q[x]$}{polynomials in finite fields Fq}}

The above construction can fail in a finite field $\FF_q[x]$ since, as already hinted, the cyclotomic polynomial $\Phi$ can become reducible. Finding an irreducible cyclotomic polynomial in $\FF_q[x]$ will to our knowledge require Artin's conjecture on primitive roots, which currently relies on the Generalised Riemann Hypothesis \cite{hooley_artins_1967}. Thus we abandon this approach and instead provide a non-constructive proof based on a counting argument.

The number of monic polynomials of degree $n$ in $\FF_q[x]$ is trivially $q^n$. On the other hand, the number of \emph{irreducible} monic polynomials $N_q(n)$ can be found in classical textbooks such as \textcite[Theorem 3.25]{lidl_finite_1996}, and is given by the formula
\begin{equation*}
    N_q(n) = \frac{1}{n} \sum_{d | n} \mu(d) q^{n/d} = \frac{1}{n} q^n + \mathcal{O}(q^{n/2})
\end{equation*}
where $\mu$ is the M\"obius function. For $n\gg q$, this expression can be very accurately estimated by its first order approximation.

Now, we pick some large $n$, and to ensure that formal derivatives will not vanish, we require that the characteristic $p$ does not divide $n$. Next, we note that
\begin{equation*}
    \frac{1}{n}q^n = q^{n - \frac{\log n}{\log q}}
\end{equation*}
and we pick $k\in\NN$ such that $N_q(n)>q^{n-k}$. In practice, $k=\ceil{\frac{\log n}{\log q}}$ works for sufficiently large $n$. We remark that $k$ grows like $\mathcal{O}(\log n)$ since $q$ is fixed.

Now, consider the generic monic polynomial of degree $n$:
\begin{equation*}
    \underbrace{a_0 + a_1 x + \dots + a_{n-k-1} x^{n-k-1}}_\text{Lower part} + a_{n-k} x^{n-k} + \dots + x^n
\end{equation*}
There are $q^{n-k}$ different lower parts among these polynomials. Meanwhile, the number of irreducible monic polynomials is strictly bigger than $q^{n-k}$ by the definition of $k$. Therefore, by the pigeonhole principle, there must exist two irreducible polynomials that share the same lower part. Call them $P$ and $Q$.

Since they share the same lower part, $P(x)-Q(x)=x^{n-k}R(x)$ for some $R$ with $\deg R = k$. Our triple will be $(Q,x^{n-k} R, P)$, which satisfies all conditions needed to be an $abc$ triple. Finally, we calculate
\begin{equation*}
    \wam(PQR,s)=
    \frac{\overbrace{1\times n^s}^\text{from $P(x)$} + \overbrace{1\times n^s}^\text{from $Q(x)$} + \overbrace{(n-k)\times1^s+\mathcal{O}(\log n)}^\text{from $x^{n-k} R(x)$}}{n^s + n^s + 1^s + \mathcal{O}(\log n) } = 
    \frac{n + 2 n^s + \mathcal{O}(\log n)}{2 n^s + 1 +\mathcal{O}(\log n)}.
\end{equation*}
The $\mathcal{O}(\log n)$ term comes from factorising $R(x)$, which we have no control over, but know it must be bounded by $\deg R = k = \mathcal{O}(\log n)$. Finally, just like earlier, we observe that $\wam$ will diverge as $n\to+\infty$ since the $n$ term in the numerator dominates all other terms in the expression.

\section[Behaviour of WAM for Re(s) greater than 1]{Behaviour of $\wam$ for $\Re(s) \ge 1$}
\label{sec:wam_analysis_numerical}

Having analysed the behaviour of $\wam(n,s)$ in the region $\Re(s) < 1$ in \cref{sec:wam_divergence_integers,sec:wam_divergence_poly}, we now turn to the region $\Re(s) \ge 1$. We provide numerical evidence based on the dataset of over $10^7$ ABC triples collected by \textcite{abcDataset}. 

\begin{figure}[ht]
    \centering
    \begin{subfigure}[b]{0.49\textwidth}
        \centering
        \includegraphics[width=\textwidth]{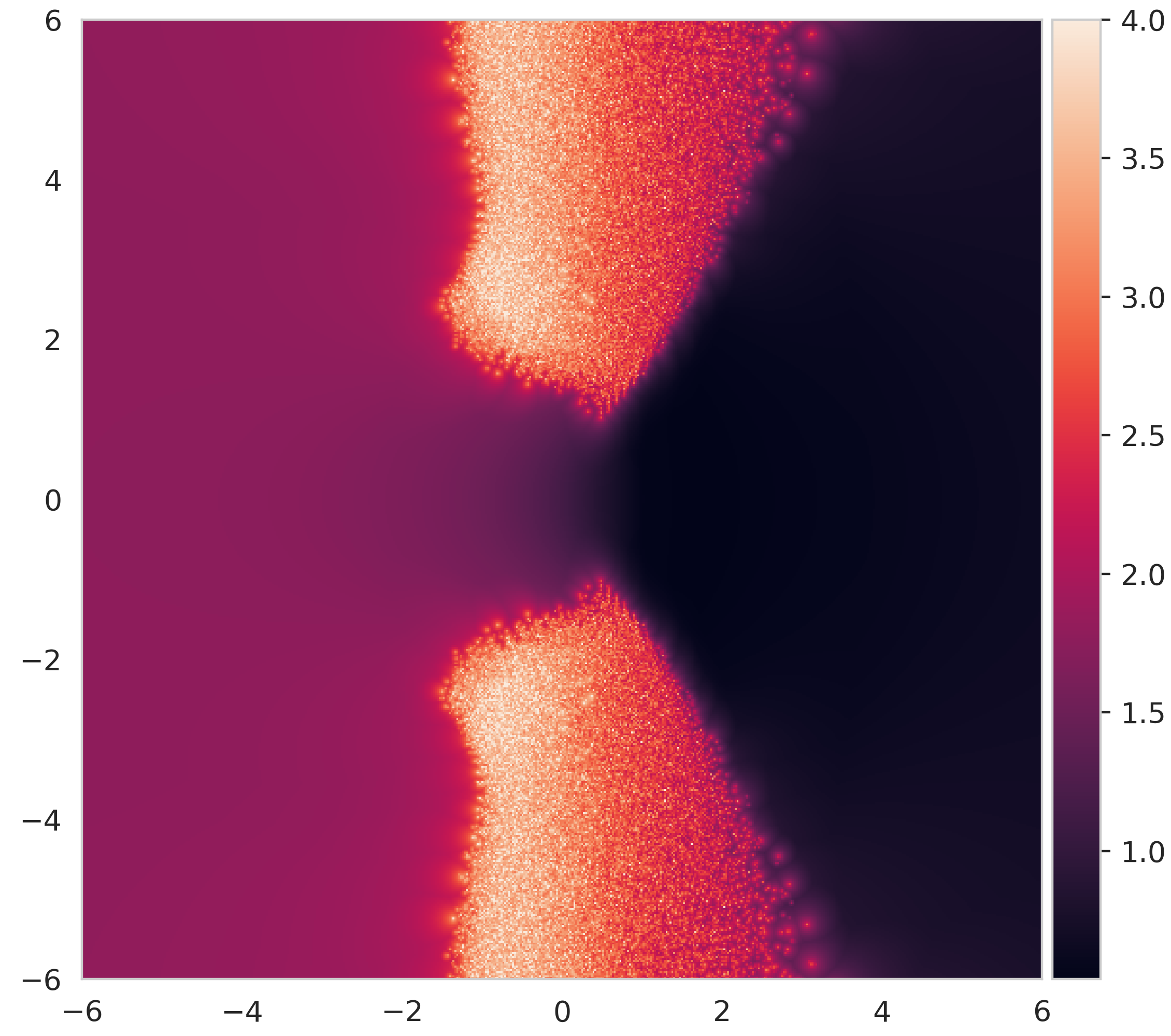}
        \caption{Large section of the complex plane, where $\max$ is taken over 50000 triples. }
        \label{fig:max_wam_complex_large}
    \end{subfigure}
    \hfill
    \begin{subfigure}[b]{0.49\textwidth}
        \centering
        \includegraphics[width=\textwidth]{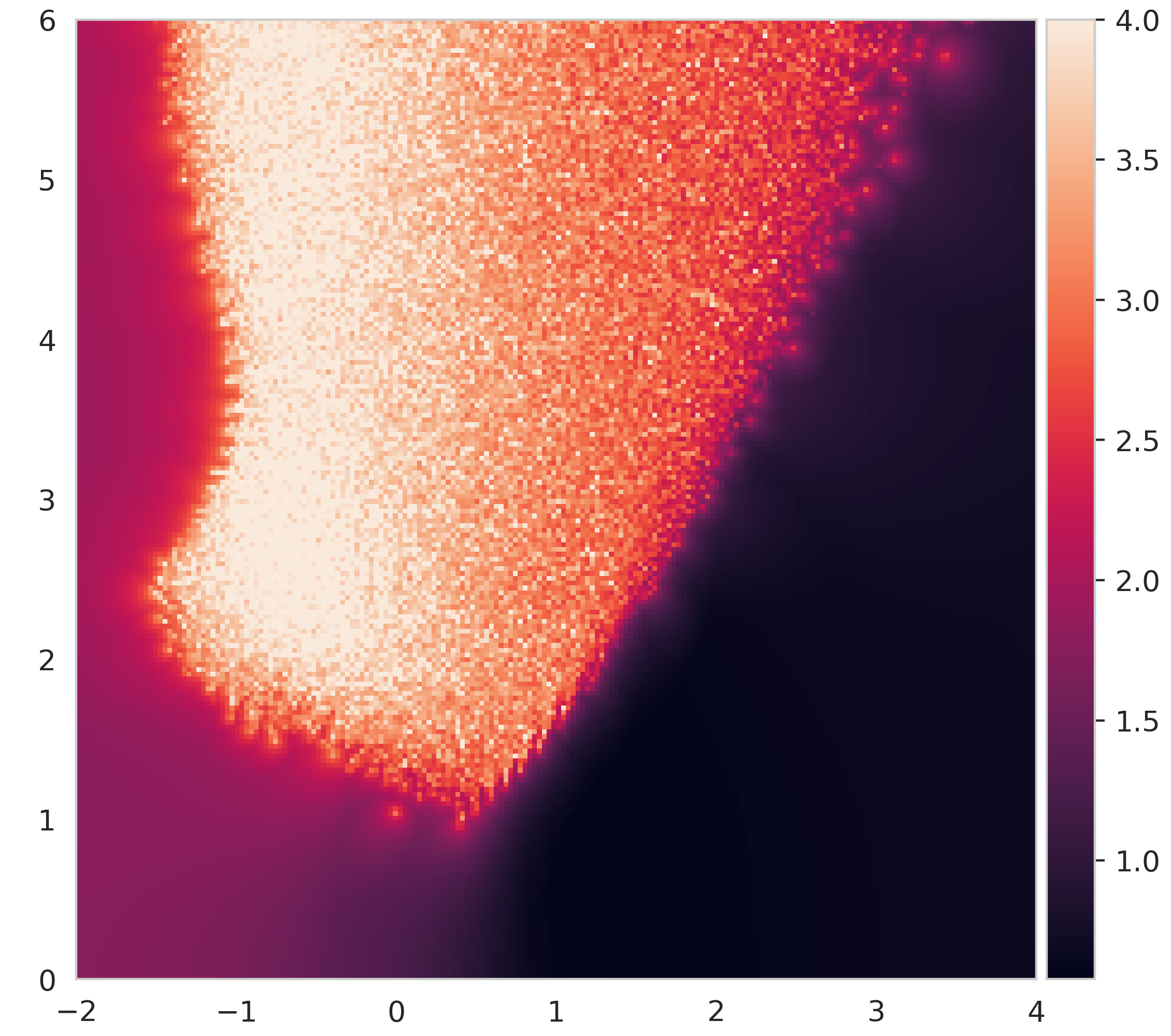}
        \caption{Zoomed in region with poles, where $\max$ is taken over $5\times 10^5$ triples.}
        \label{fig:max_wam_complex_zoom}
    \end{subfigure}
    \caption{Values of $\log_{10} \max \abs{\wam(abc,s)}$ as a function of $s$ in the complex plane, where $\max$ is taken over ABC triples with $c\sim10^{18}$. Bright colors are poles of $\wam$ for individual $n$, which are then accumulated with $\max$.}
    \label{fig:max_wam_complex_both}
\end{figure}

The quantity $\limsup|\wam(abc,s)|$ is of special interest due to its connection with the ABC conjecture, but it cannot be computed directly. Therefore, we approximate it with $\max \wam(abc,s)$, where the $\max$ is taken over the largest triples in the dataset.

\cref{fig:max_wam_complex_both} shows the behaviour of $\abs{\wam(abc,s)}$ for many large $abc$ triples, for values of $s$ whose real and imaginary parts are between $-6$ and $6$. A notable feature of this figure is the region with poles above and below the real axis. These poles have a major impact on the limiting behaviour of $\wam$, so studying them could be helpful for understanding the generalisation of the ABC conjecture to arbitrary $s$.

We will explore this further in the following sections, where we look at regions $\Re(s)=1$, $\Re(s)\to +\infty$, and $1<\Re(s)\ll +\infty$ separately.

\subsection[Behaviour for Re(s)=1]{Behaviour  for $\Re(s) = 1$}

\begin{figure}[ht]
    \begin{minipage}[t]{0.47\linewidth}
        \includegraphics[width=0.95\textwidth]{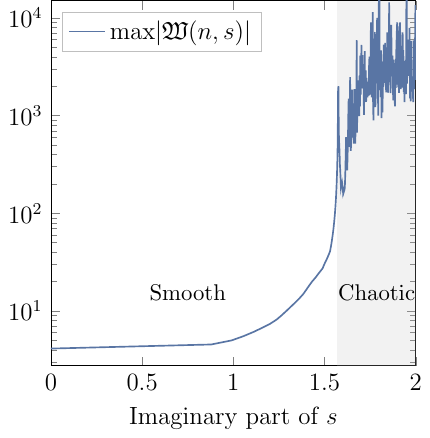}
        \caption{Values of $\wam$ on the line $\Re(s) = 1$. Chaotic behaviour occurs due to poles of $\wam$.}
        \label{fig:wam_critical_line}
    \end{minipage}
    \hfill
    \begin{minipage}[t]{0.47\linewidth}
        \includegraphics[width=\linewidth]{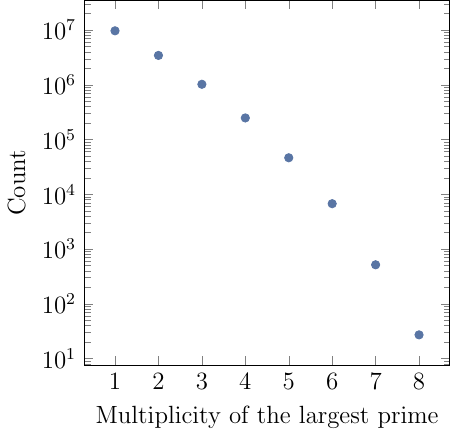}
        \caption{Number of triples with a specific multiplicity of the largest prime factor $e_m$ in the dataset.}
        \label{fig:max-prime-multiplicity}
    \end{minipage}
\end{figure}

The line $\Re(s) = 1$ has close ties to the ABC conjecture as at $s = 1$, $\wam$ simply reduces to the original WAM function. We plot $\wam$ on this line for ABC triples in \cite{abcDataset}, and show this on \cref{fig:wam_critical_line}.

The main feature of \cref{fig:wam_critical_line} is that the plot starts off with $\wam\approx3$ at $s=1$, in accordance with the ABC conjecture, but quickly explodes in magnitude and becomes chaotic. This corresponds to the region with poles in \cref{fig:max_wam_complex_zoom}.

It is natural to consider whether the smooth region shrinks in the limit, and whether it collapses entirely to $s=1$. We conjecture that this might be true, but we also note that even if the chaotic region extends to $s=1$, it does not necessarily imply divergence at the point itself.

\subsection[Behaviour for Re(s) tends to infinity]{Behaviour  for $\Re(s) \to+\infty$}

As $s$ approaches positive infinity, it is easy to see that
\begin{equation}
    \label{max_multiplicity}
    \lim_{\Re(s) \to +\infty} \wam(n,s)=e_m,
\end{equation} 
where $e_m$ is the multiplicity of the largest prime factor in $n$. The distribution of $e_m$ over $abc$ triples is unknown, so on \cref{fig:max-prime-multiplicity} we plot the distribution on the dataset of \textcite{abcDataset}.

We observe that triples with large $e_m$ are rare, but this does not exclude the possibility that $e_m$ is unbounded in the limit. 

As mentioned in the introduction, we have $e_m < \wam(n,1) \omega(n)$. If the $abc$ conjecture is true, then large $e_m$ requires large $\omega(abc)$. However, $\omega(abc)$ cannot be too large compared to $abc$ because it would imply existence of large primes in the factorisation, which would bring $e_m$ down. 


Therefore, the asymptotic behaviour of $e_m$ is closely related to, and likely as hard as, the original $abc$ conjecture.

\subsection[Behaviour for large finite Re(s)]{Behaviour for $1<\Re(s)\ll +\infty$}


Consider the primes $p_1, \dots, p_m$ appearing in factorisation of some integer $n$. Assume $n$ is not a perfect power of a squarefree number so that $\wam$ is not constant and let
\begin{equation}
    f(s) = \sum (\log p_k)^s
    \label{def:denom}
\end{equation}
be the denominator of $\wam$. Further assume that $m>2$ so that poles of $\wam$ are nontrivial. We wish to characterise the region in which $\wam(n,s)$ has no poles. A necessary condition is to find a region where $f(s)$ has no zeros. 

Let $s = a + i b$. If $\wam$ has a pole at $s$, then $f(s)=0$. Expanding $f$ one obtains 
\begin{align*}
    \sum^{m-1} (\log p_k)^s &= - (\log p_m)^s \\
    \abs*{ \sum^{m-1} (\log p_k)^s } &= \abs*{(\log p_m)^s} \\
    \sum^{m-1} (\log p_k)^a &\ge (\log p_m)^a
\end{align*}
Therefore, every pole of $\wam$ satisfies this inequality. Both left and right hand sides are monotonically increasing, but the left side is greater for small $a$ while the right side is greater for large $a$. Thus, there exists a unique $\acrit$ at which the equality is attained:
\begin{equation}
    \sum^{m-1} (\log p_k)^{\acrit} = (\log p_m)^{\acrit}
\end{equation}
and every pole $s=a+ib$ satisfies $a\le \acrit$. Thus, we have by construction the following.

\begin{prop}
    The function $f(a + i b)$ has no zeros in the region $a > \acrit$.
\end{prop}
Next, we will show that this bound is optimal and that there exist zeros arbitrarily close to the critical line $a=\acrit$. We begin by establishing some useful lemmas.

\begin{lemma}[Equidistribution]
    If $1, \alpha, \beta, ...$ are linearly independent over $\QQ$ then the set of fractional parts
    \begin{equation*}
        \braces[\bigg]{ \paren[\Big]{\braces{\alpha n}, \braces{\beta n}, \dots} \Big\vert~ n \in \NN }
    \end{equation*}
    is dense in the unit hypercube.
\end{lemma}

\begin{lemma}
    For two different primes $p$ and $q$, $\log p$ and $\log q$ are linearly independent.
    \begin{proof}
        Otherwise, one would have $p^a = q^b$ for some integers $a, b$, clearly a contradiction.
    \end{proof}
\end{lemma}

\begin{prop}
    For any set of phases $(\theta_k)_{k=1}^m$ there is a sequence $(b_t)_{t=1}^\infty$ such that for all $k$:
    \begin{equation*}
        \lim_{t\to\infty} (b_t\log p_k\ \mathrm{ mod }\ 2\pi) = \theta_k
    \end{equation*}
    \begin{proof}
        Notice that the vector $\theta/2\pi$ is inside the unit hypercube. By linear independence of $\log p_k / 2\pi$ over $\QQ$ we can find a sequence of integers $m_t$ such that for all $t$:
        \begin{equation*}
            \lim_{t\to\infty} \braces*{\frac{\log p_k}{2\pi}m_t} = \frac{\theta_k}{2\pi}
        \end{equation*}
        Then, the sequence $(b_t)_{t=1}^\infty$ defined with $b_t = 2\pi m_t$ satisfies the desired condition.
    \end{proof}
\end{prop}

\begin{corollary}
    The above sequence $(b_t)_{t=1}^\infty$ satisfies $\lim_t f(a+ib_t) = \sum (\log p_k)^a e^{i \theta_k}$.
\end{corollary}

With these statements we can now focus on the zeros of $f$ near $\acrit$.

\begin{lemma}
    There exists $\varepsilon>0$ such that the derivative of $f$ has a positive lower bound $M>0$ on the vertical strip $a\in(\acrit-\varepsilon,\acrit]$:
    \begin{equation*}
        \abs{f'(a+ib)}\ge M
    \end{equation*}
    \begin{proof}
        Let 
        \begin{equation*}
            \delta = \sum^{m-1} (\log p_k)^{a} \paren*{1 - \frac{\log\log p_k}{\log\log p_m}} > 0.
        \end{equation*}
        Next, pick $a_0$ such that 
        \begin{equation*}
            (\log p_m)^{a_0} + \delta = \sum^{m-1} (\log p_k)^{a_0}
        \end{equation*}
        This way for $a\in(a_0,\acrit]$ one has
        \begin{equation*}
            (\log p_m)^{a} + \delta > \sum^{m-1} (\log p_k)^{a}
        \end{equation*}
        Now, by computing the derivative we have 
        {\allowdisplaybreaks
        \begin{align*}
            |f'(a + ib)| &= \abs*{\sum (\log p_k)^s \log\log p_k} && \\
            &\ge \abs*{ |(\log p_m)^s| \log\log p_m - \sum^{m-1} |(\log p_k)^s| \log\log p_k } && \text{triangle ineq.} \\
            &= (\log p_m)^{a} \log\log p_m - \sum^{m-1} (\log p_k)^{a} \log\log p_k&& \\
            &> \paren*{\sum^{m-1} (\log p_k)^{a} - \delta} \log\log p_m - \sum^{m-1} (\log p_k)^{a} \log\log p_k && \text{ from } a>a_0 \\
            &= \sum^{m-1} (\log p_k)^{a} (\log\log p_m - \log\log p_k) - \delta \log\log p_m && \\
            &> \sum^{m-1} \paren[\big]{(\log p_k)^{a} -(\log p_k)^{a_0}}(\log\log p_m - \log\log p_k) && \text{ definition of }\delta \\
            &= M(a) > 0
        \end{align*}
        }
        Finally, pick some $\varepsilon<\acrit-a_0$ and $M=\min M(a)>0$.
    \end{proof}
\end{lemma}

\begin{theorem}[Rouch\'e]
    If two holomorphic functions $g$ and $h$ satisfy $|h(z)| < |g(z)|$ on some contour $\partial K$, then $g$ and $g+h$ have the same number of zeros in $K$.
\end{theorem}

\begin{prop}
\label{prop:picking_phases}
    There exists $\varepsilon>0$ such that for all $a\in(\acrit-\varepsilon,\acrit]$ and phases $(\theta_k)_{k=1}^m$ satisfying 
    \begin{equation*}
        \sum (\log p_k)^a e^{i \theta_k} = 0,
    \end{equation*}
    there is a sequence of pairs $(a_t, b_t)_{t=1}^\infty$ with:
    \begin{align*}
        \lim_{t\to\infty} a_t &= a \\ 
        \lim_{t\to\infty} (b_t\ \mathrm{ mod }\ 2\pi) &= \theta_k \mathrm{\ for\ all\ } k\\ 
        f(a_t + i b_t) &= 0 \mathrm{\ for\ all\ } t
    \end{align*}
    \begin{proof}
        By the previous proposition we have a series $(b_t)$ so that $\lim f(a+ib_t) = \sum (\log p_k)^a e^{i \theta_k} = 0$. However, $a+ib_t$ are not zeros of $f$. To fix this, we will show that there are zeros of $f$ near $a+ib_t$.

        Let $\varepsilon>0$.
        
        Let $z_0 = \acrit + i b$ for some $b$ where $|f|<\varepsilon$ and consider the two functions 
        \begin{align*}
            g(z) &= f(z_0) + f'(z_0) (z - z_0) \\
            h(z) &= f(z) - g(z)
        \end{align*}
        Since $g$ is a first order Taylor approximation to $f$, we can expect to have $|h|<|g|$ in a radius $r$ around our point. Also, $g(z)$ has a zero exactly at $z_0 - f(z_0) / f'(z_0)$. Since $|f(z_0)|<\varepsilon$ and $|f'(z_0)|>M$, we have that a zero of $g$ is within $\varepsilon/M$ of $z_0$ (which we can make arbitrarily small by shrinking $\varepsilon$ to fit inside $r$).

        Therefore, by Rouch\'e we have the $f=g+h$ also has a zero near $z_0$.

        Since these zeros will approach $a+ib_t$ arbitrarily close by shrinking $\varepsilon$, taking the sequence of these zeros completes the proof.
    \end{proof}
\end{prop}

\begin{corollary}
    The above sequence $(a_t, b_t)$, satisfies  $\lim (\log p_k)^{a_t + i b_t} = (\log p_k)^a e^{i \theta_k}$ for all $k$.
\end{corollary}

\begin{prop}
    The function $f(a + i b)$ has zeros in the $\varepsilon$-neighbourhood of the line $x=\acrit$ for all $\varepsilon>0$.
    \begin{proof}
        Set $\theta_m$ = 0 and $\theta_k=\pi$ for $k<m$. By definition of $\acrit$ one has 
        $\sum (\log p_k)^a e^{i \theta_k} = 0$. By \cref{prop:picking_phases} there exists a sequence $(a_t, b_t)$ such that $f(a_t + i b_t)=0$ and $\lim a_t = \acrit$, as required.
    \end{proof}
\end{prop}

So far we established that $f$ has zeros arbitrarily close to $a=\acrit$. However, this is not sufficient by itself to conclude that $\wam$ has poles here. For that, one has to show that at these zeroes the numerator of $\wam$ is non-zero, as the singularities could be removable otherwise.

\begin{prop}
    The function $\wam(a + i b)$ has poles in the $\varepsilon$-neighbourhood of the line $x=\acrit$ for all $\varepsilon>0$.
    \begin{proof}
        Assume for the sake of contradiction that $\acrit$ is not the optimal bound for poles of $\wam$. Then there exists $\varepsilon>0$ so that all poles with  $a \in (\acrit-\varepsilon, \acrit]$ are removable, meaning
        \begin{equation}
        \label{eq:by_contra_impl}
            f(s)=0\ \Longrightarrow\ \sum e_k (\log p_k)^s=0
        \end{equation}
        We will show that this implies all $e_k$ are the same, making $\wam$ constant and having no poles at all. To do this we will \textquote{realify} \cref{eq:by_contra_impl}, just as previously we \textquote{realified} the condition $f(s)=0$ into a statement involving $a$ only.

        Let $a_0$ be some value in range $(\acrit-\varepsilon, \acrit]$, and let $\theta_k$ be any phases that exactly satisfy the following equation:
        \begin{equation*}
            \sum (\log p_k)^{a_0} e^{i\theta_k}=0.
        \end{equation*}
        By \cref{prop:picking_phases} we can find a sequence $(a_t, b_t)$ so that $\lim a_t=a_0$, $f(a_t + ib_t)=0$ for all $t$, and each term in $f$ approaches the chosen phase individually:
        \begin{equation*}
            \lim (\log p_k)^{a_t + ib_t} = (\log p_k)^{a_0} e^{i\theta_k} \text{ for all } k
        \end{equation*}
        From \cref{eq:by_contra_impl} we know that the the numerator is zero at these points, and since it also has a well defined limit we conclude that for our choice of $a_0$ and $\theta_k$ we also have
        \begin{equation*}
            \sum e_k (\log p_k)^{a_0} e^{i\theta_k}=0.
        \end{equation*}
        By varying $a_0$ and choosing appropriate phases $\theta_k$ (namely, dependent on $a$) at each point, we obtain the following key identity:
        \begin{equation}
            \label{eq:by_contra_impl_real}
            \sum (\log p_k)^a e^{i \theta_k(a)} = 0 
            \ \Longrightarrow\ 
            \sum e_k  (\log p_k)^a e^{i \theta_k(a)} = 0
        \end{equation}
        Recall that $\sum^{m-1} (\log p_k)^a \ge (\log p_m)^a$, for all $a$ in this region, with equality holding at $\acrit$. Therefore, in the left $\varepsilon$-neighbourhood of $\acrit$ we have that $(\log p_1)^a$, $(\log p_2)^a$, and $(\log p_m)^a - \sum_3^{m-1} (\log p_k)^a$ satisfy the triangle inequality.
        

        Now, we can pick our phases $\theta_k$ so that $\theta_m=\pi$, $\theta_k=0$ for $3\le k<m$, and $\theta_1$ and $\theta_2$ are uniquely determined to satisfy the first part of \cref{eq:by_contra_impl_real}. From this we obtain:
        \begin{align}
            \label{eq:fixed_thetas_denom}
            0 &= (\log p_1)^a e^{i\theta_1(a)} + (\log p_2)^a e^{i\theta_2(a)} + \sum_3^{m-1} (\log p_k)^a - (\log p_m)^a \\
            \label{eq:fixed_thetas_num}
            0 &= e_1 (\log p_1)^a e^{i\theta_1(a)} + e_2 (\log p_2)^a e^{i\theta_2(a)} + \sum_3^{m-1} e_k (\log p_k)^a - e_m (\log p_m)^a
        \end{align}
        Multiplying \cref{eq:fixed_thetas_denom} by $e_1$ and subtracting from \cref{eq:fixed_thetas_num} we obtain
        \begin{equation*}
            0 = (e_2-e_1) \underset{\in \CC}{\underbrace{(\log p_2)^a e^{i\theta_2(a)} \vphantom{\sum_3^{m-1} (e_k-e_1) (\log p_k)^a - (e_m-e_1) (\log p_m)^a)}} 
            } + \underset{\in \RR}{\underbrace{\sum_3^{m-1} (e_k-e_1) (\log p_k)^a - (e_m-e_1) (\log p_m)^a}}
        \end{equation*}
        but this is only possible if $e_1 = e_2$ since the complex term is not constant.

        Finally, we conclude by symmetry that all $e_k$ are pairwise equal, meaning that $n$ is a $e$-th power of a squarefree number. This is a contradiction because we ruled this case out at the beginning of the section.
    \end{proof}
\end{prop}

Now that the properties of $\acrit$ have been established, we can investigate the behaviour of $\acrit$ for different $n$, and $\wam$ for $a>\acrit$.

\begin{figure}[ht]
    \centering
    \includegraphics[width=0.5\textwidth]{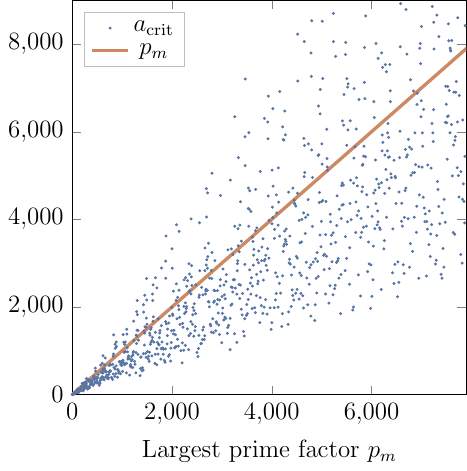}
    \caption{Numerically computed values of $\acrit$.}
    \label{fig:acrit_lim}
\end{figure}

Firstly, we numerically observe from \cref{fig:acrit_lim} that asymptotically $\acrit\sim p_m$.

Then, we note that $\wam$ is unbounded on line $a=\acrit$, even though it does not have poles on this line. Meanwhile, for any fixed $a>\acrit$ there exists a bound on $\abs\wam$ and it is derived using the triangle inequality:
\begin{equation}
\label{eq:abs_wam_bound_for_a>a_crit}
    \abs{\wam(n,a+ib)} \le \frac{\sum e_k (\log p_k)^a}{(\log p_m)^a - \sum^{m-1}(\log p_k)^a} = \wam^{upper}(a)
\end{equation}
This bound is not necessarily optimal, but it is monotonically decreasing and converges to $e_m$ in the limit, which is consistent with what was previously derived about $\wam$.

There are tools in complex analysis such as Phragm{\'e}n-Lindel{\"o}f theorem \cite{PhragLind1908} that can also produce bounds on $\wam$.

\begin{theorem}[Phragm{\'e}n-Lindel{\"o}f]
    Let $a\geq\frac{1}{2}$ and $G=\{z:|\text{arg z}|<\frac{\pi}{2a}\}$. Suppose $f$ is holomorphic on $G$ and there is a constant $M$ such that $\limsup_{z\mapsto w}{|f(z)|\leq M}$ for all $w$ in $\partial G$. If there are positive constants $P$ and $b < a$ such that
    \begin{equation*}
        |f(z)|\leq P\exp(|z|^b)
    \end{equation*}
    for all $z$ with $|z|$ sufficiently large, then $|f(z)|\leq M$ for all $z$ in $G$.
    \label{thm: PhragLind}
\end{theorem}

However, this requires establishing a bound on a boundary of some region $G$ first. By the requirement of holomorphicity, $\Re(z)$ in the region G and on the boundary $\partial G$ must be at least $a_\text{crit}$.

If the upper bound on $G$ is established using \cref{eq:abs_wam_bound_for_a>a_crit}, then the bound from Phragm{\'e}n-Lindel{\"o}f will be strictly weaker than simply applying \cref{eq:abs_wam_bound_for_a>a_crit} to the bulk of $G$.

The remaining region of interest is $1<\Re(s)< a_{\rm crit}$. We present the following collection of plots to show that there doesn't seem to be any easily discernible pattern in the poles and singularities in this region. These figures are produced by numerically computing zeros of \cref{def:denom}. In general, the behaviour seems somewhat chaotic and it is unclear to us whether more structure can be found in these patterns. 

\begin{figure}[H]
  \centering
  \begin{subfigure}[b]{0.5\textwidth}
    \includegraphics[width=\textwidth, trim={4 1 20 20}, clip]{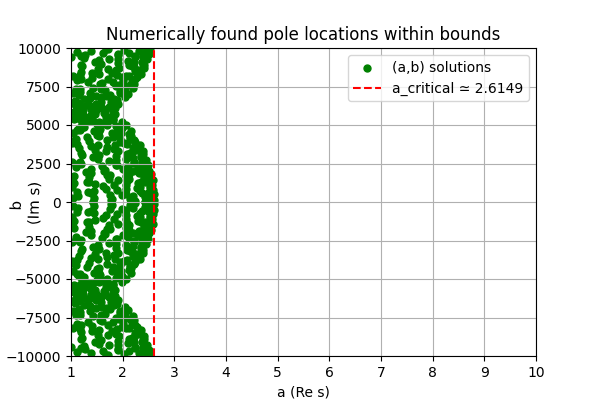}
    \caption{\parbox{\linewidth}{\centering
      $a_{\rm crit}\approx2.61, \; p_m=61$\\
      $a=1, \; b=1484734, \; c=1484735$}}
  \end{subfigure}\hfill
  \begin{subfigure}[b]{0.5\textwidth}
    \includegraphics[width=\textwidth, trim={4 1 20 20}, clip]{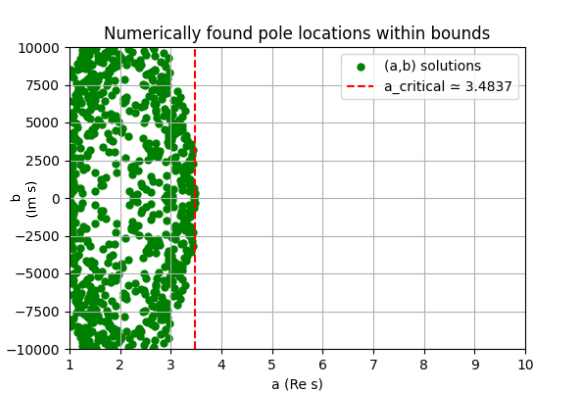}
    \caption{\parbox{\linewidth}{\centering 
      $a_{\rm crit}\approx3.48, \; p_m=19$\\
      $a=1960, \; b=59049, \; c=61009$}}
  \end{subfigure}

  \vspace{1em}

  \begin{subfigure}[b]{0.5\textwidth}
    \includegraphics[width=\textwidth, trim={4 1 20 20}, clip]{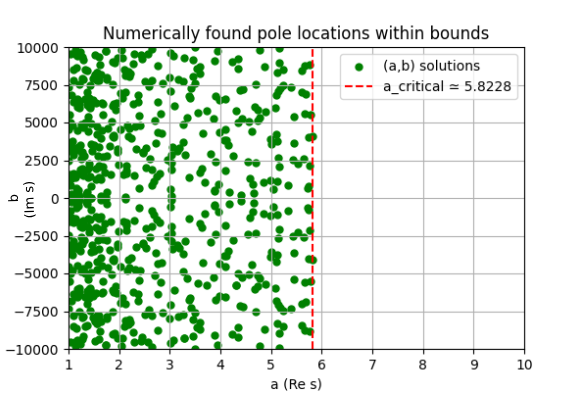}
    \caption{\parbox{\linewidth}{\centering 
      $a_{\rm crit}\approx5.82, \; p_m=29$\\
      $a=78, \; b=9765547, \; c=9765625$}}
  \end{subfigure}\hfill
  \begin{subfigure}[b]{0.5\textwidth}
    \includegraphics[width=\textwidth, trim={4 1 20 20}, clip]{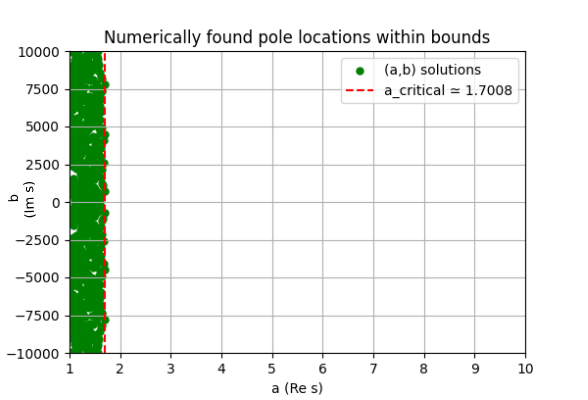}
    \caption{\parbox{\linewidth}{\centering 
      $a_{\rm crit}\approx1.70, \; p_m=521$\\
      $a=4782969, \; b=41354375, \; c=46137344$}}
  \end{subfigure}

  \vspace{1em}

  \begin{subfigure}[b]{0.5\textwidth}
    \includegraphics[width=\textwidth, trim={4 1 20 20}, clip]{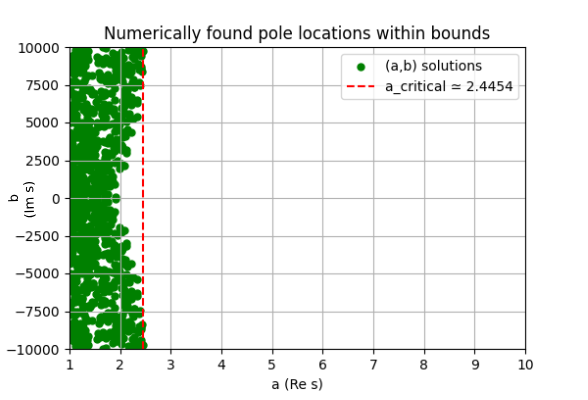}
    \caption{\parbox{\linewidth}{\centering 
      $a_{\rm crit}\approx2.45, \; p_m=113$\\
      $a=537824, \; b=134906067, \; c=135443891$}}
  \end{subfigure}\hfill
  \begin{subfigure}[b]{0.5\textwidth}
    \includegraphics[width=\textwidth, trim={4 1 20 20}, clip]{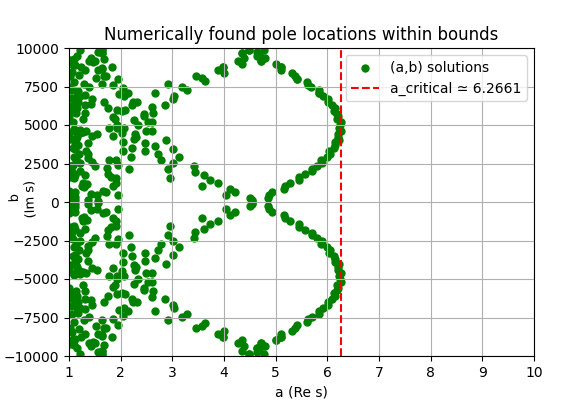}
    \caption{\parbox{\linewidth}{\centering
      $a_{\rm crit}\approx6.27, \; p_m=179$\\
      $a=13573088, \; b=349609375, \; c=363182463$}}
  \end{subfigure}

  \caption{Numerically found pole locations for six ABC triples. Each sub-figure shows the scatter plot $(a,b)$, critical vertical line at $a_{\rm crit}$ and the corresponding $abc$-triple.}
  \label{fig:pole-locations-all}
\end{figure}


\printbibliography

\end{document}